\documentclass{llncs}
        \usepackage{amsfonts}
        \usepackage{amssymb}
\usepackage{bm,bbm}
        
\newtheorem{theoreme}{Theorem}
\newtheorem{lemme}{Lemma}
\newtheorem{corollaire}{Corollary}
\def\calh{{\mathcal{H}}}
\def\calk{{\mathcal{K}}}

\def\R{{\mathbb{R}}}
\def\C{{\mathbb{C}}}
\def\E{{\mathbb{E}}}
\def\N{{\mathbb{N}}}

\def\eps{\varepsilon}
\def\be{\begin{equation}}
\def\ee{\end{equation}}

\title{Almost-Euclidean subspaces of $\ell_1^N$ via tensor products: a
simple approach to randomness reduction}
\author{Piotr Indyk\inst{1}\thanks{This research has been supported in part by David and Lucille Packard Fellowship, MADALGO (Center for Massive Data Algorithmics, funded by the Danish National Research Association) and NSF grant CCF-0728645.} \and Stanislaw Szarek\inst{2}\thanks{Supported in part by grants from the National Science Foundation (U.S.A.) and the U.S.-Israel~BSF.}}
\institute{MIT \email{indyk@mit.edu}  \and {CWRU \& Paris 6} \email{szarek@math.jussieu.fr} }

\date{}

\begin{document}

\maketitle

\begin{abstract}
    It has been known since 1970's that the $N$-dimensional $\ell_1$-space contains nearly Euclidean subspaces whose dimension is $\Omega(N)$. However, proofs of existence of such subspaces were probabilistic, hence non-constructive, which made the results not-quite-suitable for subsequently discovered applications to high-dimensional nearest neighbor search, error-correcting codes over the reals, compressive sensing and other computational problems. In this paper we present a ``low-tech'' scheme which, for any $\gamma > 0$, allows us to exhibit nearly Euclidean $\Omega(N)$-dimensional subspaces of $\ell_1^N$ while using only $N^\gamma$ random bits. Our results extend and complement (particularly) recent work by Guruswami-Lee-Wigderson. Characteristic features of our approach include (1) simplicity (we use only tensor products) and (2) yielding almost Euclidean subspaces with arbitrarily small distortions. 
\end{abstract}

\section{Introduction}

It is a well-known fact that for any vector $x \in \R^N$,  its $\ell_2$ and $\ell_1$ norms
are related by the (optimal) inequality $\|x\|_2 \le \|x\|_1 \le \sqrt{N} \|x\|_2$.
However, classical results in geometric functional analysis show that for a ``substantial fraction" of vectors , the relation between its 1-norm and 2-norm can be made much tighter.
Specifically,~\cite{FLM,Kas,GG} show that there exists a subspace $E \subset \R^N$ of dimension 
$m=\alpha N$, and a scaling constant $S$ such that for all $x \in E$
\begin{eqnarray} \label{e:l2l1}
1/D \cdot \sqrt{N} \|x\|_2 \le S \|x\|_1 \le  \sqrt{N} \|x\|_2 
\end{eqnarray}
where $\alpha \in (0,1) $ and $D=D(\alpha)$, called the {\em distortion} of $E$, are absolute (notably dimension-free) constants. Over the last few years, such ``almost-Euclidean" subspaces of $\ell_1^N$ have found numerous applications, to high-dimensional nearest neighbor search~\cite{I00}, error-correcting codes over reals and  compressive sensing~\cite{KT07,GLR08,GLW08}, vector quantization~\cite{LVer}, oblivious dimensionality reduction and $\epsilon$-samples for high-dimensional half-spaces~\cite{KRS},  and to other problems.

For the above applications, it is convenient and sometimes crucial that the subspace $E$ is defined in an explicit manner\footnote{For the purpose of this paper ``explicit" means ``the basis of $E$ can be generated by a deterministic algorithm with running time polynomial in $N$." However, the individual constructions can be even ``more explicit" than that.}. 
However, the aforementioned results do not provide much guidance in this regard, since they  use the {\em probabilistic method}. Specifically, either the vectors spanning $E$, or the vectors spanning the space dual to $E$, are i.i.d. random variables from some distribution.
As a result, the constructions require  $\Omega(N^2)$ independent random variables as starting point.
Until recently, the largest {\em explicitly} constructible almost-Euclidean subspace of $\ell_1^N$, due to 
Rudin~\cite{Rud} (cf.~\cite{LLR}),  had only a dimension  of $\Theta(\sqrt{N})$. 

During the last few years, there has been a renewed interest in the problem~\cite{AMil,Szar,I-UP,LS07,GLR08,GLW08}, with researchers using ideas gained from the study of expanders, extractors and error-correcting codes to obtain several explicit constructions.
The work  progressed on two fronts, focusing on (a) fully explicit constructions of subspaces attempting to maximize the dimension and minimize the distortion~\cite{I-UP,GLR08}, as well as (b) constructions using limited randomness, with dimension and distortion matching (at least qualitatively) the existential dimension and distortion bounds~\cite{I00,AMil,LS07,GLW08}.
The parameters of the constructions are depicted in Figure \ref{f:table}. Qualitatively, they show that in the fully explicit case, one can achieve either arbitrarily low distortion or arbitrarily high subspace dimension, but not (yet?) both.
In the low-randomness case, one can achieve arbitrarily high subspace dimension and constant distortion while using randomness that is sub-linear in $N$; achieving arbitrarily low distortion was possible as well, albeit at a price of (super)-linear randomness.

\begin{figure}
\begin{center}
\begin{tabular}{|c|c|c|c|}
\hline
Reference & Distortion & Subspace dimension & Randomness \\
\hline 
\cite{I-UP}           & $1+\epsilon$                           & $N^{1-o_{\epsilon}(1)}$ & explicit\\
\cite{GLR08}      & $(\log N)^{O_{\eta}(\log \log \log N)}$  & $(1-\eta)N$ & explicit\\
\hline
\cite{I00} &  $1+\epsilon$ & $\Omega(\epsilon^2/ \log(1/\epsilon))N $ & $O(N \log^2 N)$\\
\cite{AMil,LS07} & $O_{\eta}(1)$                          & $(1-\eta)N$ & $O(N)$ \\
\cite{GLW08}     & $2^{O_{\eta}(1/\gamma)}$                    & $(1-\eta)N$ & $O(N^{\gamma})$ \\
This paper          & $1+\epsilon$                            & $(\gamma\epsilon)^{O(1/\gamma)} N$ &  $O(N^{\gamma})$\\
\hline
\end{tabular}
\caption{The best  known results for constructing almost-Euclidean subspaces of $\ell_1^N$.
The parameters $\epsilon, \eta, \gamma \in (0,1)$ are assumed to be constants, although we explicitly point out when the dependence on them is subsumed by the big-Oh notation.}
\end{center}
\end{figure}
\label{f:table}

\paragraph{Our result} In this paper we show that, using sub-linear randomness, one can construct a subspace with  arbitrarily small distortion while keeping its dimension proportional to $N$. More precisely, we have:

\begin{theoreme}
 \label{main}
Let $\epsilon, \gamma \in (0,1)$. Given $N \in \N$, assume that we have at our disposal 
a sequence of random bits  of length $\max\{N^\gamma,C(\epsilon, \gamma)\} 
\log (N/(\epsilon\gamma))$. Then, in deterministic polynomial (in $N$) time, we can generate numbers
$M>0$, $m \geq c(\epsilon, \gamma) N$ and an $m$-dimensional subspace of 
$\ell_1^N$ $E$, for which we have
\[ \forall x \in E, \ \ \  (1-\epsilon) M \|x\|_2 \leq \|x\|_1 \leq (1+\epsilon) M \|x\|_2 \]
with probability greater than $98 \%$.
\end{theoreme}

 In  a sense, this complements the result of~\cite{GLW08}, optimizing the distortion of the subspace at the expense of its dimension.  Our approach also allows to retrieve -- 
 using a simpler and low-tech approach\ -- the results  of~\cite{GLW08} (see the comments at the end of the Introduction). 
 
 \paragraph{Overview of techniques} The ideas behind many of the prior constructions as well as this work can be viewed as variants of the related developments in the context of error-correcting codes.
  Specifically, the construction of~\cite{I-UP} resembles the approach of amplifying minimum distance of a code using expanders developed in~\cite{ABNNR}, while the constructions of~\cite{GLR08,GLW08} were inspired by  low-density parity check codes. The reason for this state of affairs is that a vector whose $\ell_1$ norms and $\ell_2$ norms are very different must be ``well-spread", i.e., a small subset of its coordinates cannot contain most of its $\ell_2$ mass (cf.~\cite{I-UP,GLR08}). This is akin to a property required from a good error-correcting code, where the weight (a.k.a. the $\ell_0$ norm) of each codeword cannot be concentrated on a small subset of its coordinates.
  
   In this vein, our construction utilizes a tool frequently used for (linear) error-correcting codes, namely the {\em tensor product}. Recall that, for two linear codes $C_1 \subset \{0,1\}^{n_1}$ and   $C_2 \subset \{0,1\}^{n_2}$, their tensor product is a code 
 $C \subset  \{0,1\}^{n_1 n_2}$, such that for any codeword $c \in C$ (viewed as an $n_1 \times n_2$ matrix), each column of $c$ belongs to $C_1$ and each row of $c$ belongs to $C_2$.
 It is known that the dimension of $C$ is a product of the dimensions of $C_1$ and $C_2$, and that the same holds for the minimum distance. 
 This enables constructing a code of ``large" block-length $N^k$ by starting from a code of ``small" block-length $N$ and tensoring it $k$ times. 
Here, we roughly show that the tensor product of two subspaces yields a subspace whose distortion is a product of the distortions of the subspaces. 
Thus, we can randomly choose an initial small low-distortion subspace, and tensor it with itself to yield the desired dimension.
 
 However, tensoring alone does not seem sufficient to give a subspace with distortion arbitrarily close to $1$. This is because we can only analyze the distortion of the product space for the case when the scaling factor $S$ in Equation~\ref{e:l2l1} is equal to $1$ (technically, we only prove the left inequality, and rely on the general relation between the $\ell_2$ and $\ell_1$ for the upper bound). For $S=1$, however, the best achievable distortion is strictly greater than $1$, and tensoring can make it only larger.  To avoid this problem, instead of the $\ell_1^N$ norm we use the $\ell_1^{N/B}(\ell_2^B)$ norm, for a ``small" value of $B$.  The latter norm (say, denoted by $\| \cdot \|$) treats the vector as a sequence of $N/B$  ``blocks" of length $B$, and returns the sum of the $\ell_2$ norms of the blocks.
We show that there exist subspaces $E \subset \ell_1^{N/B}(\ell_2^B)$ such that for any $x \in E$ we have
\begin{eqnarray*}
1/D \cdot \sqrt{N/B} \|x\|_2 \le \|x\| \le  \sqrt{N/B} \|x\|_2 
\end{eqnarray*}
for $D$ that is arbitrarily close to $1$. Thus, we can construct almost-Euclidean subspaces of $\ell_1(\ell_2)$ of desired dimensions using tensoring, and get rid of  the ``inner" $\ell_2$ norm at the end of the process.

We point out that if we do not insist on distortion arbitrarily close to $1$, the ``blocks"  
are not needed and the argument simplifies substantially. In particular, to retrieve the 
results of~\cite{GLW08},  it is enough to combine the 
scalar-valued version of Proposition \ref{mult} below with ``off-the-shelf" random 
constructions \cite{Kas,GG} yielding -- in the notation of Equation \ref{e:l2l1} -- a subspace $E$, 
for which the parameter $\alpha$ is close to $1$.

\section{Tensoring subspaces of $L_1$}

We start by defining some basic notions and notation used in this section.

\paragraph{Norms and distortion} 
In this section we adopt the ``continuous" notation for vectors and norms. 
Specifically, consider a real Hilbert space $\calh$ and a probability measure ${\mu}$ over $[0,1]$.
For $p \in [1,\infty]$ consider the space $L_p(\calh)$ of $\calh$-valued $p$-integrable functions $f$ endowed 
with the norm 
\[ \|f\|_p=\|f\|_{L_p(\calh)}= \left (\int \|f(x)\|_\calh^p \, d\mu(x)\right)^{1/p}  \]
In what follows we will omit  $\mu$ from the formulae 
since the measure will be clear from the context (and largely irrelevant).
As our main result concerns finite dimensional spaces, it suffices to focus on the case where  $\mu$ is simply the normalized counting measure over  the discrete set $\{ 0, 1/n, \ldots (n-1)/n \}$  for some fixed $n \in \mathbb{N}$
(although the statements  hold in full generality). 
%
In this setting, the functions $f$ from $L_p(\calh)$ are equivalent to $n$-dimensional vectors with coordinates in  $\calh$.\footnote{The values from $\calh$ roughly correspond to the finite-dimensional ``blocks" in the construction sketched in the introduction. Note that $\calh$  can be discretized 
similarly as the $L_p$-spaces; alternatively, functions that are constant on intervals of the type $\big((k-1)/N, k/N\big)$ can be considered in lieu of discrete measures.} 
The advantage of using the $L_p$ norms as opposed to the $\ell_p$ norms that the relation between the $1$-norm and the $2$-norm does not involve scaling factors that depend on dimension, i.e., we have  $\|f\|_2 \geq \|f\|_1$ for all $f \in  L_2(\calh)$ (note that, for the  $L_p$ norms,  the ``trivial'' inequality goes in the other direction than for the $\ell_p$ norms).   This simplifies the notation considerably.

We will be interested in lialmost
 subspaces $E \subset L_2(\calh)$ on which the $1$-norm and $2$-norm uniformly agree, i.e., 
for some $c \in (0,1]$,
\be \label{21est}
\|f\|_2 \geq \|f\|_1 \geq c \|f\|_2
\ee
for all $f \in E$. The best (the largest) constant $c$ that works in (\ref{21est}) 
will be denoted $\Lambda_1(E)$. For completeness, we also define 
$\Lambda_1(E)=0$ if no $c>0$ works.

\paragraph{Tensor products} If $\calh, \calk$ are Hilbert spaces, 
$\calh\otimes_2\calk$ is their 
(Hilbertian) tensor product, which may be (for example) described by the 
following property: if $(e_j)$ is an orthonormal sequence in $\calh$ and $(f_k)$ is an orthonormal sequence in $\calk$, then  $(e_j\otimes f_k)$ is an orthonormal sequence in 
$\calh\otimes_2\calk$ (a basis if $(e_j)$ and $(f_k)$ were bases). 
Next, any element of $L_2(\calh)\otimes \calk$ is canonically identified 
with a function in the space $L_2(\calh\otimes_2 \calk)$; 
note that such functions are $\calh\otimes\calk$-valued, but are defined 
on the same probability space as their counterparts from $L_2(\calh)$. If $E \subset L_2(\calh)$ is 
a linear subspace, $E\otimes \calk$ is -- under this identification -- 
a linear subspace of $L_2(\calh\otimes_2 \calk)$.

\medskip As hinted in the Introduction, our argument depends (roughly) on the fact that 
the property  expressed by (\ref{e:l2l1}) or (\ref{21est}) ``passes" 
to tensor products of subspaces, and that it ``survives" replacing scalar-valued 
functions by ones that have values in a Hilbert space.  Statements to similar effect 
of various degrees of generality and precision are widely available in the 
mathematical literature, see for example \cite{MZ,Bec,And,FJ}. However, 
we are not aware of a reference that subsumes all the facts needed 
here and so we present an elementary self-contained proof.

\medskip We start with two preliminary lemmas.
\begin{lemme} \label{lemme1} 
If $g_1, g_2, \ldots \in E \subset  L_2(\calh)$, then 
$$
 \int \big( \sum_k \|g_k(x)\|_\calh^2\big)^{1/2} \, dx \geq  
 \Lambda_1(E) \, \Big( \int  \sum_k \|g_k(x)\|_\calh^2\, dx \Big)^{1/2}.
 $$
\end{lemme}
{\em Proof }\ Let  $\calk $ be an auxiliary Hilbert space and  $(e_k)$ an 
orthonormal sequence (O.N.S.) in $\calk$. 
We will apply Minkowski inequality -- a continuous version of the triangle inequality, 
which says that for vector valued functions $\|\int h \| \leq \int \|h\|$ -- 
to the $\calk$-valued function $h(x) = \sum_k \|g_k(x)\|_\calh \, e_k$. 
As is easily seen, $\|\int h\|_{\calk} = \|\sum_k \big(\int \|g_k(x)\|_\calh \, dx\big)\, e_k\|_{\calk}
= \big(\sum_k \|g_k\|_{L_1(\calh)} ^2\big)^{1/2}$. 
Given that $g_k \in E$, $\|g_k\|_{L_1(\calh)}  \geq 
\Lambda_1(E)\, \|g_k\|_{L_2(\calh)} $ and so 
$$
\Big\|\int h\Big\|_{\calk} \geq  
\Lambda_1(E) \Big( \int  \sum_k \|g_k(x)\|_\calh^2\, dx \Big)^{1/2}$$
On the other hand, the left hand side of the inequality in Lemma \ref{lemme1} is 
exactly $ \int \|h\|_{\calk} $, so the Minkowski inequality yields the required 
estimate.

\medskip We are now ready to state the next lemma. Recall that $E$ is a linear subspace of  $L_2(\calh)$, and  $\calk $ is a Hilbert space.
\begin{lemme} \label{lemme2}
$\Lambda_1(E\otimes \calk) = \Lambda_1(E)$
\end{lemme}
If $E \subset L_2=L_2(\mathbb{R})$, the lemma says that any estimate of type  
(\ref{21est}) for scalar functions $f \in E$ carries over to their 
linear combinations with vector coefficients, 
namely to functions of the type $\sum_j v_j f_j$, $f_j \in E, v_j\in \calk$. 
In the general case, any estimate for $\calh$-valued functions 
$f \in E \subset L_2(\calh)$ carries over to functions of the form 
$\sum_j  f_j\otimes v_j \in L_2(\calh\otimes_2 \calk)$, with $f_j \in E, v_j\in \calk$. 

\medskip \noindent {\em Proof of Lemma \ref{lemme2}} \  
Let $(e_k)$ be an orthonormal basis  of $\calk$. 
In fact w.l.o.g. we may assume that $\calk = \ell_2$ and that $(e_k)$ 
is the canonical  orthonormal basis.
Consider $g=\sum_j  f_j\otimes v_j$, where $f_j \in E$ and $v_j\in \calk$. 
Then also $g=\sum_k  g_k\otimes e_k$  for some $g_k \in E$ and 
hence (pointwise) 
$\|g(x)\|_{\calh\otimes_2\calk} =\big( \sum_k \|g_k(x)\|_\calh^2\big)^{1/2}$.
Accordingly, 
$\|g\|_{L_2(\calh\otimes_2 \calk)}= \big( \int  \sum_k \|g_k(x)\|_\calh^2\, dx \big)^{1/2}$, 
while 
$\|g\|_{L_1(\calh\otimes_2 \calk)}=  \int \big( \sum_k \|g_k(x)\|_\calh^2\big)^{1/2} \, dx $. 
Comparing such quantities is exactly the object of Lemma \ref{lemme1}, 
which implies that 
$\|g\|_{L_1(\calh\otimes_2 \calk)} \geq  
 \Lambda_1(E)\|g\|_{L_2(\calh\otimes_2 \calk)}.$ 
 Since $g \in E\otimes \calk$ was arbitrary, it follows that 
 $\Lambda_1(E\otimes \calk) \geq \Lambda_1(E)$. 
 The reverse inequality is automatic (except in the trivial case 
 $\dim \calk =0$, which we will ignore).

\medskip If $E \subset L_2(\calh)$ and $F \subset L_2(\calk)$ are subspaces, 
$E\otimes F$ is the subspace of $L_2(\calh\otimes_2 \calk)$ spanned by 
$f\otimes g$ with $f\in E, g\in F$. (For clarity, $f\otimes g$ is a function 
on the {\em product } of the underlying probability spaces and is defined by 
$(x,y) \to f(x)\otimes g(y) \in \calh\otimes\calk$.) 

\medskip
 The next proposition shows the key property of tensoring almost-Euclidean spaces.

\begin{proposition} \label{mult}
$\Lambda_1(E\otimes F) \geq \Lambda_1(E) \Lambda_1(F)$
\end{proposition}
{\em Proof }\ 
Let $( \varphi_j )$ and $( \psi_k )$ be orthonormal bases of respectively $E$ and $F$ 
and let $g = \sum_{j,k} t_{jk}  \, \varphi_j \otimes \psi_k$. 
We need to show that 
$\|g\|_{L_1(\calh\otimes_2 \calk)} 
\geq \Lambda_1(E) \Lambda_1(F) \|g\|_{L_2(\calh\otimes_2 \calk)}$, 
where the $p$-norms refer to the product probability space, for example
$$
\|g\|_{L_1(\calh\otimes_2 \calk)} 
= \int \int 
\big\|\sum_{j,k} t_{jk}  \, \varphi_j(x) \otimes \psi_k(y)\big\|_{\calh\otimes_2 \calk}\; dx \, dy .
$$
Rewriting the expression under the sum and subsequently applying 
Lemma \ref{lemme2} to the inner integral for fixed $y$ gives
\begin{eqnarray*}
\int \big\|\sum_{j,k} t_{jk}  \, \varphi_j(x) \otimes \psi_k(y)\big\|_{\calh\otimes_2 \calk}\; dx 
&=& \int \big\|\sum_{j}\varphi_j(x)  \otimes \Big(\sum_{k}t_{jk}  \,  
\psi_k(y)\Big)\big\|_{\calh\otimes_2 \calk}\; dx \\
&\geq&\Lambda_1(E) \Big(\int \big\|\sum_{j}\varphi_j(x)  
\otimes \Big(\sum_{k}t_{jk}  \,  \psi_k(y)\Big)\big\|_{\calh\otimes_2 \calk}^2\; dx \Big)^{1/2}\\
&=&\Lambda_1(E) \Big(\sum_{j}  \big\|\sum_{k}t_{jk}  \,  \psi_k(y)\big\|_{ \calk}^2 \Big)^{1/2}
\end{eqnarray*}
In turn, $\sum_{k}t_{jk}  \,  \psi_k \in F$ (for all $j$) and so, by Lemma \ref{lemme1},
\begin{eqnarray*}
\int \Big(\sum_{j}  \big\|\sum_{k}t_{jk}  \,  \psi_k(y)\big\|_{ \calk}^2 \Big)^{1/2}\, dy 
&\geq& \Lambda_1(F)
\Big( \int \sum_{j}  \big\|\sum_{k}t_{jk}  \,  \psi_k(y)\big\|_{ \calk}^2 \, dy\Big)^{1/2} \\
&=&  \Lambda_1(F)\;  \|g\|_{L_2(\calh\otimes_2 \calk)} .
\end{eqnarray*} 
Combining the above formulae yields the conclusion of the Proposition.

\medskip

\section{The construction}

In this section we describe our low-randomness construction. We start from a recap of the probabilistic construction, since we use it as a building block.

\subsection{Dvoretzky's theorem, and its ``tangible" version}

For general normed spaces, the following is one possible statement of the well-known Dvoretzky's theorem \cite{dv}: 

\medskip \noindent {\em Given $m\in \N$ and $\eps > 0$ there is $N=N(m,\eps)$ such that, for any norm on $\R^N$ there is an $m$-dimensional subspace on which the ratio of $\ell_1$ and $\ell_2$ norms is (approximately) constant, up to a multiplicative factor $1+\eps$.} 

\medskip \noindent For specific norms this statement can be made more precise, both in describing the dependence 
$N=N(m,\eps)$ and in identifying the constant of (approximate) proportionality of norms. The following version is (essentially) due to Milman \cite{milman}. 

\medskip \noindent 
{\bf Dvoretzky's theorem } (Tangible version) \ {\sl Consider the $N$-dimensional Euclidean space (real or complex) endowed with the Euclidean norm $\| \cdot \|_2$ and some other norm  $\| \cdot \|$ such that, for some $b>0$, $ \|\cdot\| \leq b \|\cdot\|_2$. 
Let $M =  \E\|X\|$, where $X$ is a random variable uniformly distributed on the unit Euclidean sphere. 
Then there exists a computable universal constant $c >0$, so that if $0<\eps <1$ and $m \leq c\eps^2 (M/b)^2 N$, then for more than 99\% (with respect to the Haar measure) $m$-dimensional subspaces $E$ 
we have
\begin{equation} \label{dvor}
\forall x \in E, \ \ \  (1-\eps) M \|x\|_2 \leq \|x\| \leq (1+\eps) M \|x\|_2 .
\end{equation}
}
Alternative good expositions of the theorem are in, e.g.,  \cite{FLM}, \cite{MS} and \cite{pisier}.
We point out that standard and most elementary proofs yield $m \leq c\eps^2 /\log(1/\eps) (M/b)^2 N$; the dependence on $\eps$ of order $\eps^2$ was obtained in the important papers \cite{Gor,Sch}, see also \cite{ASW}. 

\subsection{The case of $\ell_1^n(\ell_2^B)$}

Our objective now is to apply Dvoretzky's theorem and subsequently Proposition \ref{mult} to spaces of the form $\ell_1^n(\ell_2^B)$ for some $n, B \in \N$, so from now on we set $\|\cdot\| := \| \cdot\|_{\ell_1^n(\ell_2^B)}$
To that end, we need to determine the values of the parameter $M$ that appears in the theorem. (The optimal value of $b$ is clearly $\sqrt{n}$, as in the scalar case, i.e., when $B=1$.) We have the following standard (cf. \cite{Ball}, Lecture 9)

\begin{lemme} \label{ave}
\[
M(n,B):= \E_{x \in S^{nB -1}}\ \|x\| = 
\frac{\Gamma(\frac{B+1}2)}{\Gamma(\frac{B}2)}
\frac{\Gamma(\frac{nB}2)}{\Gamma(\frac{nB+1}2)} \ n .
\]
In particular,  $\sqrt{ 1 + \frac 1{n-1}} \sqrt{\frac 2\pi} \, \sqrt{n} > M(n,1) > \sqrt{\frac 2\pi} \, \sqrt{n}$ for all $n \in \N$ (the scalar case) and 
$M(n,B) > \sqrt{1 - \frac 1B} \, \sqrt{n}$ for all $n, B \in \N$.
\end{lemme}

The equality is shown by relating (via passing to polar coordinates) spherical averages of norms to Gaussian means: if $X$ is a random variable uniformly distributed on the Euclidean sphere $S^{N-1}$ and $Y$ has the standard  Gaussian distribution on $\R^N$, then, for any norm 
$\|\cdot\|$, 
\[
\E\|Y\|\ = \ \frac{\sqrt{2}\, \Gamma(\frac{N+1}2)}{\Gamma(\frac{N}2)}\; \E\|X\|
\]
The inequalities follow from the estimates $\sqrt{x - \frac 12} < \frac {\Gamma(x + \frac 12)}{\Gamma(x)} < \sqrt{x}$ (for $x \geq {\frac 12}$), which in turn are consequences of 
log-convexity of $\Gamma$ and its functional equation 
$\Gamma(y+1)=y\Gamma(y)$.  (Alternatively, Stirling's formula may be used to arrive at a similar conclusion.)

\medskip Combining Dvoretzky's theorem with Lemma \ref{ave}  yields

\begin{corollaire} \label{onestep:cont}
If $0<\eps <1$ and $m \leq c_1 
\eps^2  n$, then for more than $99\%$  of the  $m$-dimensional subspaces $E \subset \ell_1^n$ we have
\begin{equation}
\label{e:c1}
 \forall x \in E \ \ \  (1-\eps) \sqrt{\frac 2 \pi}\,\sqrt{n}  \|x\|_2 \leq \|x\|_1 \leq  (1+\eps) \sqrt{ 1 + \frac 1{n-1}}\sqrt{\frac 2 \pi}\sqrt{n} \|x\|_2 
 \end{equation}
Similarly, if $B>1$ and $m \leq 
c_2 \eps^2 nB$, then  for  more than $99\%$  of the $m$-dimensional subspaces $E \subset \ell_1^n(\ell_2^B)$ we have
\begin{equation}
\label{e:c2}
\forall x \in E \ \ \  (1-\eps) \sqrt{1-\frac 1B}\,\sqrt{n}  \|x\|_2 \leq \|x\| \leq \sqrt{n} \|x\|_2 
\end{equation}

\end{corollaire}

We point out that the upper estimate on $\|x\|$ in the second inequality is valid for all $x \in \ell_1^n(\ell_2^B)$ and, like the estimate $M(n,B) \leq \sqrt{n}$, 
follows just from the Cauchy-Schwarz inequality.

\medskip
Since  a random subspace chosen uniformly according to the Haar measure on the manifold of $m$-dimensional subspaces of $\R^N$ (or $\C^N$) 
can be constructed from an $N \times m$ random Gaussian matrix, we may apply 
 standard discretization techniques to obtain the following
\begin{corollaire} \label{onestep}
There is a deterministic algorithm that, given  $\eps, B, m,n$ as in Corollary~\ref{onestep:cont} and a sequence of $O(m n \log (mn/\epsilon))$ random bits, generates subspaces $E$ as in Corollary~\ref{onestep:cont} with probability greater than $98 \%$, in time polynomial in $1/\eps+B+m+n$.
\end{corollaire}

We point out that in the literature on the ``randomness-reduction", one typically uses
Bernoulli matrices  in lieu of Gaussian ones.
This enables avoiding the discretization issue, since the problem is phrased directly in terms of random bits.  
Still, since proofs of Dvoretzky type theorems for Bernoulli matrices 
are often much harder than for their Gaussian counterparts, we prefer 
to appeal instead to a simple discretization of Gaussian random 
variables.  We note, however, that the early  approach of \cite{Kas} 
was based on Bernoulli matrices.

\medskip
We are now ready to conclude the proof of Theorem \ref{main}.  
Given $\eps \in (0,1)$ and $n \in \N$, choose $B = \lceil \eps^{-1}\rceil$ and $m= \lfloor c\eps^2 (1-\frac 1B) nB\rfloor \geq c_0 \eps^2 nB$. 
Corollary \ref{onestep} (Equation~\ref{e:c2}) and repeated application of Proposition \ref{mult}  give us a subspace 
$F \subset \ell_1^\nu(\ell_2^\beta)$ (where $\nu=n^k$ and $\beta = B^k$) of dimension 
$m^k \geq (c_0  \eps^2)^k\nu \beta$ such that 
\[ \forall x \in F \ \ \  (1-\eps)^{3k/2} n^{k/2}  \|x\|_2 \leq \|x\| \leq  n^{k/2}  \|x\|_2 .\]
Moreover, $F=E\otimes E \otimes \ldots \otimes E$, where $E \subset \ell_1^n(\ell_2^B)$ is a typical $m$-dimensional subspace. Thus
in order to produce $E$, hence $F$, we only need to generate a 
``typical" $m \approx c_0  \eps^2 (\nu \beta))^{1/k}$  subspace 
of the $nB=(\nu \beta))^{1/k}$-dimensional space $\ell_1^n(\ell_2^B)$.  
Note that for fixed $\eps$ and $k>1$, $nB$ and $m$ are asymptotically 
(substantially) smaller than $\dim F$.  
Further, in order to efficiently represent $F$ as a subspace of an 
$\ell_1$-space, we only need to find a good embedding of $\ell_2^\beta$ 
into $\ell_1$. This can be done using Corollary \ref{onestep} (Equation~\ref{e:c1});
note that $\beta$ depends only on $\eps$ and $k$.  
Thus we reduced the problem of finding ``large" almost Euclidean subspaces 
of $\ell_1^N$ to similar problems for much smaller dimensions. 

Theorem~\ref{main} now follows from the above discussion.
The argument gives, e.g., $c(\eps, \gamma) = (c\eps \gamma)^{3/\gamma}$ and  
$C(\eps, \gamma)=c(\eps, \gamma)^{-1}$.

\newcommand{\etalchar}[1]{$^{#1}$}
\def\cprime{$'$} \def\cprime{$'$}

\end{document}